\pdfoutput=1
\RequirePackage{ifpdf}
\ifpdf 
\documentclass[pdftex]{sigma}
\else
\documentclass{sigma}
\fi

\usepackage[all,arc]{xy}
\usepackage{multirow}
\usepackage{yhmath}
\numberwithin{equation}{section}

\newtheorem{Theorem}{Theorem}[section]
\newtheorem{Proposition}[Theorem]{Proposition}
\newtheorem{Conjecture}[Theorem]{Conjecture}
\newtheorem{Question}[Theorem]{Question}

{ \theoremstyle{definition}
	\newtheorem{Definition}[Theorem]{Definition}}

\begin{document}

\allowdisplaybreaks

\newcommand{\arXivNumber}{2006.00745}

\renewcommand{\thefootnote}{}

\renewcommand{\PaperNumber}{078}

\FirstPageHeading

\ShortArticleName{The Fundamental Groups of Open Manifolds with Nonnegative Ricci Curvature}

\ArticleName{The Fundamental Groups of Open Manifolds\\ with Nonnegative Ricci Curvature\footnote{This paper is a~contribution to the Special Issue on Scalar and Ricci Curvature in honor of Misha Gromov on his 75th Birthday. The full collection is available at \href{https://www.emis.de/journals/SIGMA/Gromov.html}{https://www.emis.de/journals/SIGMA/Gromov.html}}}

\Author{Jiayin PAN}

\AuthorNameForHeading{J.~Pan}

\Address{Department of Mathematics, University of California-Santa Barbara,\\ Santa Barbara CA 93106, USA}
\Email{\href{mailto:j_pan@math.ucsb.edu}{j\_pan@math.ucsb.edu}}
\URLaddress{\url{http://web.math.ucsb.edu/~j_pan/}}

\ArticleDates{Received June 02, 2020, in final form August 04, 2020; Published online August 17, 2020}

\Abstract{We survey the results on fundamental groups of open manifolds with nonnegative Ricci curvature. We also present some open questions on this topic.}


\Classification{53C21; 53C23; 57S30}

\renewcommand{\thefootnote}{\arabic{footnote}}
\setcounter{footnote}{0}

\section{Introduction}\label{sec_intro}

One fascinating topic in Riemannian geometry is the interplay between curvature and topology. Our primary interest in this survey is how nonnegative Ricci curvature determines the structure of fundamental groups of open (complete and non-compact) manifolds.

For comparison, we recall the topological implications of nonnegative sectional curvature for open manifolds. By Cheeger--Gromoll soul theorem \cite{CG_soul}, any open manifold $M$ of $\operatorname{sec}\ge 0$ is homotopic to a closed totally geodesic submanifold in $M$. In particular, $M$ has finite topology. On the fundamental group level, the soul theorem implies that $\pi_1(M)$ is finitely generated and virtually abelian (contains an abelian subgroup of finite index). In contrast, an open manifold with $\operatorname{Ric}\ge 0$ may have infinite second or higher Betti number, even when the manifold has Euclidean volume growth \cite{AKL,Menguy,Per_ex,ShaYang1,ShaYang2,Wraith1,Wraith2}; the reader may refer to Section 4 of \cite{SS} for a survey on the examples of manifolds with nonnegative Ricci curvature and infinite topology. Therefore, it is particularly interesting to see whether the structure of fundamental groups, finite generation and virtual abelianness, still holds for open manifolds with nonnegative Ricci curvature.

On finite generation, this is the longstanding Milnor conjecture raised in 1968.

\begin{Conjecture}[Milnor \cite{Mil}]\label{Milnor_conj}
Let $M$ be an open $n$-manifold of $\operatorname{Ric}\ge 0$. Then $\pi_1(M)$ is finitely generated.
\end{Conjecture}

The Milnor conjecture remains open.

On virtual abelianness, Wei constructed an open manifold of positive Ricci curvature whose fundamental group is not virtually abelian \cite{Wei}. Later, Wilking showed that any finitely generated virtually nilpotent group can be realized the fundamental group of some open manifold of positive Ricci curvature~\cite{Wilk}.
On the other hand, by the work of Gromov and Milnor~\cite{Gro2, Mil}, nonnegative Ricci curvature implies that any finitely generated subgroup of the fundamental group is virtually nilpotent. Kapovitch and Wilking proved a stronger result by dropping the finite generation condition and achieving a uniform bound on the index~\cite{KW}.

We summarize these results in the table below.
\begin{table}[h]
	\centering
	\begin{tabular}{|c|c|c|c|}
		\hline
		\multicolumn{2}{|c|}{let $M$ be an open manifold of} & $\operatorname{sec}\ge 0$ & $\operatorname{Ric}\ge 0$ \\
		\hline
		\multirow{3}*{then $\pi_1(M)$ is} & finitely generated & true & unknown\\
		\cline{2-4}
	 & virtually abelian & true & false\\
		\cline{2-4}
		 & virtually nilpotent & true & true\\
		\hline
	\end{tabular}\par
	\caption{Nonnegative curvature and fundamental groups.}
\end{table}

Given that proving or disproving the Milnor conjecture seems very difficult with the current understanding of nonnegative Ricci curvature, it is natural to ask the following question:

\begin{Question}\label{main_ques_1}
	For an open manifold of $\operatorname{Ric}\ge 0$, on what conditions is $\pi_1(M)$ finitely generated?
\end{Question}

Though virtual abelianness does not hold for all open manifolds with nonnegative Ricci curvature, we can ask a similar question:

\begin{Question}\label{main_ques_2}
	For an open manifold of $\operatorname{Ric}\ge 0$, on what conditions is $\pi_1(M)$ virtually abelian?
\end{Question}

In Section \ref{sec_basic}, we go through some basic tools in Ricci curvature and their applications to fundamental groups. This covers Bishop--Gromov relative volume comparison, Cheeger--Gromoll splitting theorem, and structure results on asymptotic cones. For applications on fundamental groups, we include the finiteness result by Anderson~\cite{An} and Li~\cite{Li}, Milnor's estimate on the growth function~\cite{Mil}, and fundamental groups of closed manifolds with nonnegative Ricci curvature~\cite{CG_split}.

We survey results related the Milnor conjecture in Section \ref{sec_finite_gen}. We describe Gromov's short ge\-ne\-ra\-tors~\cite{Gro1}, Kapovitch--Wilking's bound on the number of short generators \cite{KW}, and Wilking's reduction on the Milnor conjecture~\cite{Wilk}. Regarding answers to Question \ref{main_ques_1}, we survey the work by Sormani on manifolds with small linear diameter growth \cite{Sor}, the classification by Liu in dimension~$3$~\cite{Liu}, and the recent progress by the author via asymptotic geometry \cite{Pan2,Pan1,Pan3}.

In Section \ref{sec_abel}, we first describe the construction by Wei~\cite{Wei} and other related examples by Nabonnand \cite{Nab} and B\'erard-Bergery~\cite{Bergery}. Then we present author's recent work on Question~\ref{main_ques_2} \cite{Pan3,Pan4,Pan5}.

We also post some further questions and conjectures in Sections \ref{sec_finite_gen} and \ref{sec_abel}.

We mention that there are also many beautiful results on the finite topology of open manifolds under additional geometric constraints like sectional curvature, conjugate radius, diameter growth, etc. We are unable to survey these results here due to the limited space.

\section{Some basic tools in Ricci curvature}\label{sec_basic}

We review some basic tools in Ricci curvature with a focus on the applications to fundamental groups of open manifolds with $\operatorname{Ric}\ge 0$. Though some of the results in this section hold in a~more general setting, we will only consider the case of nonnegative Ricci curvature.

Because the fundamental group acts isometrically on the Riemannian universal cover by path lifting, these tools are commonly applied to not only the manifold itself but also its universal cover when we study the fundamental group.

\subsection{Volume comparison}\label{sec_vol_comp}

One of the basic tools is the Bishop--Gromov relative volume comparison.
\begin{Theorem}[Bishop--Gromov \cite{Bishop,Gro_book}]\label{volume_comp}
	Let $M$ be an open $n$-manifold of $\operatorname{Ric}\ge 0$ and let $x\in M$. Then
\[\mathrm{vol}(B_r(x))\le \mathrm{vol}(B^n_r(0)),\qquad \dfrac{\mathrm{vol}(B_R(x))}{\mathrm{vol}(B_r(x))}\le \dfrac{\mathrm{vol}(B^n_R(0))}{\mathrm{vol}(B^n_r(0))}=\dfrac{R^n}{r^n}\]
	for all $R>r>0$, where $B_r(x)$ is the metric $r$-ball centered at $x$ and $B_r^n(0)$ is the $r$-ball in the standard Euclidean space $\mathbb{R}^n$.
\end{Theorem}

Theorem \ref{volume_comp} implies that on an open manifold $M$ of $\operatorname{Ric}\ge 0$, the limit
\[\lim\limits_{r\to\infty}\dfrac{\mathrm{vol}(B_r(x))}{\mathrm{vol}(B_r^n(0))}=L\]
always exists with $L\in[0,1]$, where $x\in M$. Moreover, this limit does not depend on the choice of $x$. When $L>0$, we say that $M$ has \textit{Euclidean volume growth of constant}~$L$.

If $L$ equals the maximum value $1$, then $M$ is isometric to the standard Euclidean space~$\mathbb{R}^n$. When $M$ has almost max volume growth, that is, $L$ is sufficiently close to $1$, there is a~corresponding almost rigidity result: Perelman proved that such a manifold is contractible~\cite{Per}; later Cheeger and Colding proved that it is diffeomorphic to $\mathbb{R}^n$ \cite{CC1}.

For $L>0$, Li and Anderson independently proved that $\pi_1(M)$ is finite.

\begin{Theorem}[Li \cite{Li}, Anderson \cite{An}]\label{Eu_vol}
	Let $M$ be an open $n$-manifold of $\operatorname{Ric}\ge 0$. If $M$ has Euclidean volume growth of constant $L$, then $\pi_1(M)$ is finite with order at most $1/L$.
\end{Theorem}

In particular, if $L>1/2$, then $M$ is simply connected.
Li's proof uses heat kernel on the universal cover. We briefly go through Anderson's proof, which utilizes the volume comparison and a method to calculate the volume via Dirichlet domain. As $\pi_1(M,x)$ acts freely, discretely, and isometrically on the Riemannian universal cover $\widetilde{M}$, we can construct the Dirichilet domain~$D$ centered at $\tilde{x}$, where $\tilde{x}\in\widetilde{M}$ is a lift of~$x$. $D$ is a fundamental domain; in particular, under $\pi_1(M,x)$-action, $\gamma\cdot(B_r(\tilde{x})\cap D)$ and $\gamma'\cdot(B_r(\tilde{x})\cap D)$ are disjoint whenever $\gamma$ and $\gamma'$ are different elements in $\pi_1(M,x)$. Also, $D$ satisfies $\mathrm{vol}(B_r(x))=\mathrm{vol}(B_r(\tilde{x})\cap D)$ for all $r>0$. To prove Theorem \ref{Eu_vol}, we take a finite subset $S\subseteq \pi_1(M,x)$ and choose $d>0$ such that $d(\gamma\tilde{x},\tilde{x})\le d$ for all $\gamma\in S$. Let
\[ S(k)=\{\gamma\in\pi_1(M,x) \,|\,\gamma \text{ can be expressed as a word in } S \text{ with word length}\le k\},\]
where $k\in \mathbb{N}$. Note that $d(\gamma\tilde{x},\tilde{x})\le kd$ for all $\gamma\in S(k)$. For any $r>0$, we estimate the volume of $S(k)\cdot (B_r(\tilde{x})\cap D)$:
\begin{gather*}
\#S(k)\cdot \mathrm{vol}(B_r(x))=\#S(k)\cdot \mathrm{vol}(B_r(\tilde{x})\cap D)=\mathrm{vol}(S(k)\cdot (B_r(\tilde{x})\cap D))
\le \mathrm{vol}(B_{r+kd}(\tilde{x})).
\end{gather*}
It follows from volume comparison that
\[ \#S(k)\le \dfrac{\mathrm{vol}(B_{r+kd}(\tilde{x}))}{\mathrm{vol}(B_r(x))}\le\dfrac{\mathrm{vol}(B^n_{r+kd}(0))}{\mathrm{vol}(B_r(x))}.
\]
Let $r\to\infty$, we see that
\[ \#S(k)\le\lim\limits_{r\to\infty} \dfrac{\mathrm{vol}(B^n_{r+kd}(0))}{\mathrm{vol}(B_r^n(0))}\cdot\dfrac{\mathrm{vol}(B_r^n(0))}{\mathrm{vol}(B_r(x))}=\dfrac{1}{L}.
\]
Since this estimate holds for any $k\in\mathbb{N}$ and any finite subset $S\subseteq \pi_1(M,x)$, Theorem~\ref{Eu_vol} follows. This proof also indicates that if $M$ and $\widetilde{M}$ have the same volume growth rate, then $\pi_1(M)$ is finite as well.

Later, this finiteness result was extended to any manifold whose volume growth has order $n-\epsilon$ for some sufficiently small $\epsilon>0$.

\begin{Theorem}[Wu \cite{Wu}]\label{almost_eu_vol}
Given $n$, there exists a constant $\epsilon(n)>0$ such that the following holds.
	
Let $M$ be an open $n$-manifold of $\operatorname{Ric}\ge 0$. If there is $x\in M$ such that
\[\lim\limits_{r\to\infty} \dfrac{\mathrm{vol}(B_r(x))}{r^{n-\epsilon(n)}}\in (0,\infty), \]
	then $\pi_1(M)$ is finite.
\end{Theorem}

Also see Wan's work on finiteness under a relative volume growth condition \cite{Wan}, which implies Theorem~\ref{almost_eu_vol}. It is unknown whether the volume growth condition in Theorem~\ref{almost_eu_vol} can be weakened to $\liminf$ or any $\epsilon\in(0,1)$.

When $\pi_1(M,x)$ is infinite, Milnor proved that the volume comparison provides an estimate on the growth function of any finitely generated subgroup, which in turn implies virtual nilpotency by Gromov's work.

\begin{Theorem}\label{vir_nilp}
	Let $M$ be a complete $($open or closed$)$ $n$-manifold of $\operatorname{Ric}\ge 0$. Then any finitely generated subgroup of $\pi_1(M)$
\begin{enumerate}\itemsep=0pt
\item[$(1)$] has polynomial growth with degree at most $n$ $($Milnor {\rm \cite{Mil})},
\item[$(2)$] is virtually nilpotent $($Gromov {\rm \cite{Gro2})}.
\end{enumerate}
\end{Theorem}

We can derive Theorem \ref{vir_nilp}(1) from Anderson's proof. Let $S$ be a generating set of the finitely generated subgroup. By taking $r=1$, then we can estimate $\#S(k)$ by
\[\#S(k)\le \dfrac{\mathrm{vol}(B^n_{1+kd}(0))}{\mathrm{vol}(B_1(x))},\]
where the right hand side is a polynomial of $k$ with degree $n$.

Theorem \ref{vir_nilp}(2) is useful to prove that certain fundamental group is virtually abelian under additional conditions: without lose of generality, we can start with a nilpotent one. It also follows from Theorem~\ref{vir_nilp}(2) that every finitely generated subgroup is finitely presented~\cite{Wylie}.

We end this section with Kapovitch--Wilking's result, which significantly improves Theo\-rem~\ref{vir_nilp}(2).

\begin{Theorem}[Kapovitch--Wilking \cite{KW}]\label{C_nilp}
Given $n$, there exists a constant $C(n)$ such that the following holds.
	
Let $M$ be an open $n$-manifold of $\operatorname{Ric}\ge 0$. Then $\pi_1(M)$ contains a nilpotent subgroup $N$ of index at most $C(n)$. Moreover, any finitely generated subgroup of $N$ has nilpotency length at most $n$.
\end{Theorem}

\subsection{Splitting theorem}\label{sec_split}

A line is an isometric embedding $\gamma\colon \mathbb{R}\to M$; in other words,
\[ d(\gamma(s),\gamma(t))=|s-t|\]
for all $s,t\in\mathbb{R}$. Cheeger--Gromoll splitting theorem is a fundamental tool in nonnegative Ricci curvature.

\begin{Theorem}[Cheeger--Gromoll \cite{CG_split}]\label{CG_split}
	Let $M$ be an open $n$-manifold of $\operatorname{Ric}\ge 0$. If $M$ contains a line, then $M$ splits isometrically as a metric product $\mathbb{R}\times N$.
\end{Theorem}

The splitting theorem implies the virtual abelianness of $\pi_1(M)$ when $M$ is closed and has nonnegative Ricci curvature.

\begin{Theorem}[Cheeger--Gromoll \cite{CG_split}]\label{cpt_vir_abel}
Let $M$ be a closed $n$-manifold of $\operatorname{Ric}\ge 0$. Then $\pi_1(M)$ is virtually abelian.
\end{Theorem}

In fact, as a consequence of Theorem \ref{CG_split}, the universal cover of $M$ must split isometrically as $\mathbb{R}^k\times N$, where $N$ is a compact manifold. Consequently, its isometry group splits as well
\[ \operatorname{Isom}\big(\widetilde{M}\big)=\operatorname{Isom}\big(\mathbb{R}^k\big)\times \operatorname{Isom}(N),\]
where $\operatorname{Isom}(N)$ is a compact Lie group. Also, without lose of generality, we can assume that~$\pi_1(M)$ is nilpotent by Theorem~\ref{vir_nilp}. Since we can view $\pi_1(M)$ as a subgroup of $\operatorname{Isom}\big(\widetilde{M}\big)$, Theorem~\ref{cpt_vir_abel} follows from the algebraic fact that any nilpotent subgroup of $\operatorname{Isom}\big(\mathbb{R}^k\big)\times K$ is virtually abelian, where $K$ is a compact Lie group.

We mention that it was conjectured by Fukaya and Yamaguchi that the index of the abelian subgroup in Theorem~\ref{cpt_vir_abel} can be uniformly bounded by some constant $C(n)$~\cite{FY}. This conjecture is open even for nonnegative sectional curvature.

Theorem~\ref{cpt_vir_abel} has an analog for open manifolds. For each element $\gamma\in\pi_1(M,x)$, among all the loops based at $x$ representing $\gamma$, we can choose one of the minimal length. This loop is a~geodesic loop and we call it a \textit{representing geodesic loop} of~$\gamma$.

\begin{Theorem}\label{bdd_vir_abel}
	Let $M$ be an open $n$-manifold of $\operatorname{Ric}\ge 0$. If for some $x\in M$, all representing geodesic loops of elements in $\pi_1(M,x)$ are contained in a bounded set, then $\pi_1(M)$ is virtually abelian.
\end{Theorem}

Though in this case, the universal cover may not be the metric product of $\mathbb{R}^k$ and a compact space, it can be shown that $\pi_1(M,x)$ is contained in $\operatorname{Isom}\big(\mathbb{R}^k\big)\times K$ for some compact Lie group~$K$, which is sufficient for virtual abelianness (see Appendix~A of~\cite{Pan4} for a proof).

The assumption in Theorem \ref{bdd_vir_abel} always holds when $\operatorname{sec}\ge 0$ by setting $x$ in a soul of $M$. This follows from Cheeger--Gromoll soul theorem and Sharafutdinov retraction \cite{CG_soul,Shara}. On the other hand, the assumption in Theorem \ref{bdd_vir_abel} is very restricted for manifolds with nonnegative Ricci curvature. In fact, if $M$ has positive Ricci curvature and an infinite fundamental group, then it follows from the Cheeger--Gromoll splitting theorem that the representing geodesic loops $c_\gamma$ will always escape from any bounded sets as $\gamma$ exhausts $\pi_1(M,x)$~\cite{SW}. In Section~\ref{sec_escape}, we will see the author's work \cite{Pan4,Pan5} generalizing Theorem~\ref{bdd_vir_abel}.

\subsection{Asymptotic geometry}\label{sec_asym}

The Gromov--Hausdorff distance between two metric spaces measures how close they look alike. Using relative volume comparison Theorem \ref{volume_comp}, Gromov proved a precompactness result. This concept of convergence revolutionizes the Riemannian geometry.

\begin{Theorem}[Gromov \cite{Gro_book}]\label{precpt}
Let $(M_i,x_i)$ be a sequence of complete $n$-manifolds of $\operatorname{Ric}\ge -(n-1)$, then after passing to a subsequence, it converges in the $($pointed$)$ Gromov--Hausdorff topology to a length metric space $(X,x)$.
\end{Theorem}

The limit space $X$ is called a \textit{Ricci limit space}. Cheeger, Colding, and Naber developed the fundamental theory on the structure of these Ricci limit spaces \cite{CC1,CC2,CC3,CC4,ChN1,ChN2,Co,CoN}.

In the context of an open manifold $M$ with $\operatorname{Ric}\ge 0$, we can use Gromov--Hausdorff convergence to study its asymptotic geometry. For any sequence $r_i\to\infty$, passing to a subsequence if necessary, we obtain a convergent sequence
\[\big(r_i^{-1}M,x\big)\overset{{\rm GH}}\longrightarrow (Y,y).\]
The limit $(Y,y)$ is called an \textit{asymptotic cone of}~$M$, or a \textit{tangent cone of~$M$ at infinity}. It does not depend on the base point $x$, but it may depend on the scaling sequence~$r_i$.

Cheeger and Colding substantially generalized the splitting theorem to Ricci limit spaces~\cite{CC1}. For convenience, we denote $\mathcal{M}(n,0)$ the set of all Ricci limit spaces coming from sequences of complete $n$-manifolds of $\operatorname{Ric}\ge 0$.

\begin{Theorem}[Cheeger--Colding \cite{CC1}]\label{CC_split}
	Let $X\in\mathcal{M}(n,0)$. If $X$ contains a line, then $X$ splits isometrically as a metric product $\mathbb{R}\times N$.
\end{Theorem}

One important class of asymptotic cones is metric cones.

\begin{Theorem}[Cheeger--Colding \cite{CC1}]\label{metric_cone}
Let $M$ be an open $n$-manifold of $\operatorname{Ric}\ge 0$. Suppose that~$M$ has Euclidean volume growth. Then any asymptotic cone $(Y,y)$ of~$M$ is a metric cone with vertex~$y$.
\end{Theorem}

For a metric cone $C(Z)\in\mathcal{M}(n,0)$, $\operatorname{diam}(Z)\ge \pi$ is equivalent to $C(Z)$ containing a line. Then by Theorem~\ref{CC_split}, $C(Z)$ splits isometrically as $\mathbb{R}^k\times C(Z')$, where $\operatorname{diam}(Z')<\pi$.

Since we want to study $\pi_1(M,x)$, which acts on the universal cover $\widetilde{M}$ as isometries, it is natural to apply the asymptotic method to $\big(\widetilde{M},\pi_1(M,x)\big)$ as well. This involves the equivariant Gromov--Hausdorff convergence, introduced by Fukaya and Yamaguchi \cite{Fky,FY}. Let $r_i\to\infty$, after passing to a subsequence, we can obtain the diagram below:
\vspace{3pt}
\begin{gather}
\begin{CD}
	\big(r^{-1}_i\widetilde{M},\tilde{x},\pi_1(M,x)\big) @>{\rm GH}>>
	\big(\widetilde{Y},\tilde{y},G\big)\\
	@VV\pi V @VV\pi V\\
	\big(r^{-1}_iM,x\big) @>{\rm GH}>> \big(Y=\widetilde{Y}/G,y\big).
	\end{CD}\label{star}
\end{gather}
The $\pi_1(M,x)$-actions on $r_i^{-1}\widetilde{M}$ passes a limit $G$-action on $\widetilde{Y}$, where $G$ is some closed subgroup of $\operatorname{Isom}\big(\widetilde{Y}\big)$. It follows from Theorem~\ref{isom_Lie} below that $G$ is a Lie group.

\begin{Theorem}[Colding--Naber \cite{CoN}]\label{isom_Lie}
	Let $X\in\mathcal{M}(n,0)$. Then $\operatorname{Isom}(X)$ is a Lie group.
\end{Theorem}

Theorem \ref{isom_Lie} was first proved for the non-collapsing case by Cheeger--Colding~\cite{CC3}, then the general case by Colding--Naber~\cite{CoN}.

If $X\in\mathcal{M}(n,0)$ is a metric cone, then due to the splitting $X=\mathbb{R}^k\times C(Z')$, where $Z'$ has diameter less than $\pi$, the isometry group of $X$ splits as well
\[ \operatorname{Isom}(X)=\operatorname{Isom}\big(\mathbb{R}^k\big)\times \operatorname{Isom}(Z').\]

\section{Finite generation}\label{sec_finite_gen}

\subsection{Gromov's short generators and their properties}\label{sec_short_gen}

Gromov introduced a geometric method to choose generators of $\pi_1(M,x)$ successively \cite{Gro1}. We first choose $\gamma_1\in \pi_1(M,x)-\{e\}$ such that
\[d(\gamma_1\tilde{x},\tilde{x})=\min_{\gamma\in \pi_1(M,x)-\{e\}} d(\gamma\tilde{x},\tilde{x}).\]
Supposing that $\gamma_k$ is chosen for some $k\in\mathbb{N}$, we pick $\gamma_{k+1}\in \pi_1(M,x)-H_k$ such that
\[d(\gamma_{k+1}\tilde{x},\tilde{x})=\min_{\gamma\in \pi_1(M,x)-H_k} d(\gamma\tilde{x},\tilde{x}),\]
where $H_k=\langle\gamma_1,\dots ,\gamma_k \rangle$ is the subgroup generated by first $k$ generators. It follows from the choice of these generators that
\[d(\gamma_i\tilde{x},\gamma_j\tilde{x})\ge \max\{d(\gamma_i\tilde{x},\tilde{x}),d(\gamma_j\tilde{x},\tilde{x}) \}\]
for all $i$, $j$. When $M$ has $\operatorname{sec}\ge 0$, applying Toponogov comparison theorem to the geodesic triangle with vertices $\tilde{x}$, $\gamma_i\tilde{x}$, and $\gamma_j\tilde{x}$, one can see that $\theta_{i,j}$, the angle of this geodesic triangle at~$\tilde{x}$, is at least $\pi/3$. A standard packing argument then leads to Gromov's bound on the number of short generators.

\begin{Theorem}[Gromov \cite{Gro1}]\label{sec_sg}
	Given $n$, there exists a constant $C(n)$ such that the following holds.
	
	Let $M$ be an open $n$-manifold of $\operatorname{sec}\ge 0$ and let $x\in M$. Then the number of short generators of $\pi_1(M,x)$ is bounded by $C(n)$.
\end{Theorem}

For Ricci curvature, Kapovitch and Wilking proved the following:

\begin{Theorem}[Kapovitch--Wilking \cite{KW}]\label{ric_sg}
Given $n$, there exists a constant $C(n)$ such that the following holds.
	
Let $M$ be an open $n$-manifold of $\operatorname{Ric}\ge 0$ and let $x\in M$. Let $\Gamma$ be a finitely generated subgroup of $\pi_1(M,x)$. Then there is a point $q\in B_1(x)$ such that the number of short generators of~$\Gamma$ at~$q$ is bounded by~$C(n)$.
\end{Theorem}

The statement proved by Kapovitch and Wilking is indeed local. Note that unlike Theorem~\ref{sec_sg}, Theorem~\ref{ric_sg} gives a bound on the number of short generators at a point near~$x$, not exactly at $x$. The proof of Theorem~\ref{ric_sg} is much more complicated than the sectional curvature case. It starts with a contradicting sequence
\[(M_i,x_i)\overset{{\rm GH}}\longrightarrow(X,x)\]
so that the number of short generators at all points in $B_1(x_i)$ is larger than $i$. The proof involves an induction argument on the dimension of $X$ in the Colding--Naber sense and a successive blow-up argument to increase the dimension of the limit space. In this process, one has to shift base points multiple times. Thus this method cannot bound the number of short generators at a~priorly fixed reference point.

\looseness=-1 Note that if a bound $C(n)$ can be obtained at the base point $x$, then finite generation would follow directly by taking the first $C(n)+1$ many short generators to generate a subgroup. Also see the work by Rong and the author where a bound at the base point is achieved under certain conditions on $\pi_1(M,x)$-action \cite{PR}. On the other hand, Theorem \ref{ric_sg} does immediately imply that any finitely generated subgroup of $\pi_1(M)$ can be generated by at most $C(n)$ many elements. However, this is far away from the Milnor conjecture. A typical example with this algebraic property is the set of rational numbers $\mathbb{Q}$ as an additive group: $\mathbb{Q}$ is not finitely generated, but every finitely generated subgroup is cyclic. An even simpler example is the cyclic dyadic rationals
\[\Gamma=\left\{\dfrac{p}{2^i}\, \big|\, p\in\mathbb{Z},\,i\in\mathbb{N} \right\}.\]
A set of short generators of $\Gamma$ could potentially be $\gamma_1=1/2$, \dots, $\gamma_i=1/\big(2^i\big)$, etc. Ruling out these abelian groups as fundamental groups is crucial. In fact, with the help of Wilking's reduction, it is sufficient to consider abelian ones.

\begin{Theorem}[Wilking \cite{Wilk}]\label{reduction}
Let $M$ be an open $n$-manifold of $\operatorname{Ric}\ge 0$. If $\pi_1(M)$ is not finitely generated, then $\pi_1(M)$ has an abelian subgroup that is not finitely generated.
\end{Theorem}

Another geometric property of short generators is about its representing geodesic loops. Let $c$ be a unit speed representing geodesic loop of a short generator $\gamma$ and let $r$ be the length of $c$, then by the choice of short generators, $c$ must be minimal up to its halfway, that is, $d(x,m)=r/2$, where $m=c(r/2)$ is the midpoint of $c$. This implies that $m$ is a cut point of $x$; in particular, for $y \in M$ with
$d(x,y) > r/2$, we have
\[d(x,m)+d(y,m)-d(x,y)>0.\]
Applying the Abresch--Gromoll inequality \cite{AG} to the universal cover, Sormani established a quantitative version of the above estimate~\cite{Sor}:
for all $y\in M$ with $d(x,y)\ge r/2+S(n)r$, we have
\[d(x,m)+d(y,m)-d(x,y)\ge 2S(n)r,\]
where $S(n)>0$ is a small constant depending only on $n$. If $\pi_1(M,x)$ is not finitely generated, then $M$ has the above \textit{uniform cut} estimate at a sequence of midpoints $m_i$ with $d(x,m_i)\to\infty$. With this, Sormani confirmed the Milnor conjecture when $M$ has small linear diameter growth.

\begin{Theorem}[Sormani \cite{Sor}]\label{small_diam}
	Given $n$, there is a constant $S(n)>0$ such that the following holds.
	
	Let $M$ be an open $n$-manifold of $\operatorname{Ric}\ge 0$.
	If $M$ has small linear diameter growth
\[\limsup_{r\to\infty} \dfrac{\operatorname{diam}(\partial B_r(x))}{r}\le S(n),\]
	then $\pi_1(M)$ is finitely generated.
\end{Theorem}

The diameter here is extrinsic so $\operatorname{diam}(\partial B_r(x))\le 2r$ always holds. Theorem~\ref{small_diam} also covers manifolds with linear volume growth \cite{Sor2}, that is,
\[\limsup_{r\to\infty} \dfrac{\mathrm{vol}(B_r(x))}{r}<\infty.\]

If $\pi_1(M,x)$ is not finitely generated, then we have a sequence of scales $r_i=d(\gamma_i\tilde{x},\tilde{x})\to\infty$, where $\{\gamma_i\}_i$ is a set of short generators. Using this sequence $r_i$, we can see the consequences of infinite generation through asymptotic geometry.

\begin{Proposition}\label{sg_infty}
	Let $M$ be an open $n$-manifold of $\operatorname{Ric}\ge 0$. Suppose that $\pi_1(M)$ is not finitely generated with a set of short generators $\{\gamma_1,\dots ,\gamma_i,\dots \}$. Then in the diagram~\eqref{star} with $r_i=d(\gamma_i\tilde{x},\tilde{x})\to\infty$,
\begin{enumerate}\itemsep=0pt
\item[$(1)$] $Y$ is not polar at $y$ {\rm \cite{Sor}};
\item[$(2)$] the orbit $G\cdot \tilde{y}$ is not connected {\rm \cite{Pan1}}.
\end{enumerate}
\end{Proposition}

$Y$ being polar at $y$ means that for any $z\in Y-\{y\}$, there is a ray starting at $y$ going through $z$. Recall that short generators $\{\gamma_1,\dots ,\gamma_i,\dots \}$ provides a sequence of representing geodesic loops based at $x$, on which the uniform cut estimate holds at the midpoint~$m_i$. Passing this estimate to the asymptotic cone $Y$, we get a limit cut point in $Y$. This leads to~(1) above. For~(2), let~$H_i$ be the subgroup of $\pi_1(M,x)$ that is generated by the first $i$ many short generators. Then we consider the convergence
\[ \big(r_i^{-1}\widetilde{M},\tilde{x},H_i,\pi_1(M,x)\big )\overset{{\rm GH}}\longrightarrow \big(\widetilde{Y},\tilde{y},H,G\big).\]
The orbit $H\cdot \tilde{y}$ contains the connected component of $G\cdot \tilde{y}$ that includes~$\tilde{y}$. $\gamma_{i+1}$ converges to an element $\gamma_\infty\in G$ with $d(H\cdot \tilde{y},\gamma_\infty \tilde{y})=1$; in particular, $G\cdot \tilde{y}$ is not connected.

We mention that Sormani observed another way to see the consequence of infinite generation from afar by using the \textit{loop to infinity} property~\cite{Sor3}. An open manifold satisfies the loop to infinity property, if for any noncontractible loop $c$ and any compact subset $K$, there exists a~loop in $M-K$ that is freely homotopic to~$c$. This property gives a surjective inclusion map $i_\star\colon \pi_1(M-K,y)\to \pi_1(M,y)$, where $K$ is compact and $y\in M-K$. Sormani proved that for any open manifold~$M$ of $\operatorname{Ric}\ge 0$, either~$M$ has the loop to infinity property or a double cover of~$M$ splits off a~line isometrically; in particular, manifolds of positive Ricci curvature always have this property. As a result, if there is an example~$M$ with $\operatorname{Ric}>0$ and an infinitely generated~$\pi_1(M)$, then we can slide all these infinitely many generators to the end of~$M$.

\subsection{3-manifolds}\label{sec_3}

The Milnor conjecture is true in dimension $3$. Schoen and Yau proved that any open $3$-manifold of positive Ricci curvature is diffeomorphic to $\mathbb{R}^3$ \cite{SchYau}. Using Schoen--Yau's minimal surface theory and Perelman's solution to Poincar\'e conjecture, Liu classified open $3$-manifold of nonnegative Ricci curvature \cite{Liu}: either it is diffeomorphic to $\mathbb{R}^3$ or its universal cover splits off a line isometrically; in particular, this confirms the Milnor conjecture in dimension $3$.

\begin{Theorem}[Liu \cite{Liu}]\label{3_mfds}
The Milnor conjecture holds in dimension $3$.
\end{Theorem}

Later, the author gave an alternative proof of the finite generation of $\pi_1(M)$ by using asymptotic geometry and short generators~\cite{Pan1}. Here we give an outline of the proof by the author. With Wilking's reduction and a topological result in dimension~$3$ by Evans and Moser~\cite{EM}, if~$\pi_1(M)$ is not finitely generated, then without lose generality we can assume that it is a~subgroup of~$\mathbb{Q}$; in particular, $\pi_1(M)$~is torsion-free. This implies a geometric feature on the asymptotic cones: for any equivariant asymptotic cone $ \big(\widetilde{Y},\tilde{y},G\big)$ of $\big(\widetilde{M},\Gamma\big)$,
the limit orbit $G\cdot y$ cannot be discrete.

Suppose that $\pi_1(M,x)$ is not finitely generated. Let $\{\gamma_1,\dots ,\gamma_i,\dots \}$ be a set of short generators and let $r_i=d(\gamma_i\tilde{x},\tilde{x})\to\infty$. We consider the corresponding asymptotic cones as in the diagram~\eqref{star}. To derive a contradiction, the proof eliminates all the possibilities of asymptotic cones $\widetilde{Y}$ and~$Y$ regarding their dimensions in the Colding--Naber sense~\cite{CoN}, which are integers~$1$,~$2$ or~$3$.

\textit{Case 1.} $\dim\big(\widetilde{Y}\big)=3$. This case occurs if and only if $\widetilde{M}$ has Euclidean volume growth. Hence~$\widetilde{Y}$ is a metric cone. If~$\widetilde{Y}$ splits off an $\mathbb{R}^2$ or $\mathbb{R}^3$ factor, then Cheeger--Colding theory implies that $\widetilde{M}$ is isometric to~$\mathbb{R}^3$~\cite{CC2}. If~$\widetilde{Y}$ splits off exactly one~$\mathbb{R}$ factor, then the orbit $G\cdot\tilde{y}$ is contained in a line. It follows from Theorem~\ref{sg_infty}(2) that the orbit must be discrete, which contradicts the non-discreteness that we derived from the reduction. If~$\widetilde{Y}$ does not split off any line, then the orbit $G\cdot \tilde{y}$ must be the cone tip $\widetilde{y}$; a contradiction to both non-discreteness from the reduction and the non-connectedness from Theorem~\ref{sg_infty}(2).

\textit{Case 2.} $\dim\big(\widetilde{Y}\big)=\dim(Y)=2$. Since $Y=\widetilde{Y}/G$, it is reasonable to believe that $G$ is a discrete group, thus the orbit $G\cdot y$ must be discrete as well. This again contradicts the non-discreteness. Note that ruling out this case does not require that $Y$ and $\widetilde{Y}$ come from the scaling sequence $r_i$, the length of short generators.

\textit{Case 3.} $\dim(Y)=1$. A Ricci limit space of dimension $1$ is either a ray or a line \cite{Honda}. Combined with Proposition~\ref{sg_infty}(1), the $(Y,y)$ must be a ray with $y$ not being the starting point. By choosing a suitable rescaling $s_i\to\infty$, we can obtain a different asymptotic cone $Y'$ with $\dim(Y')\ge 2$. Together with the results from cases 1 and 2, one can derive a contradiction as well.

\subsection{Equivariant stability at infinity}\label{sec_equi_stab}

With the idea of investigating equivariant asymptotic cones, the author proved new partial results on the Milnor conjecture \cite{Pan2,Pan3}. These results are related to certain geometric stability condition on the universal cover at infinity and the question below:

\begin{Question}\label{equi_stab_infty}
	\textit{Does geometric stability at infinity implies equivariant stability at infinity?}
\end{Question}

To understand how Question \ref{equi_stab_infty} is related to the Milnor conjecture, we think of an example: the universal cover has a unique asymptotic cone as a standard Euclidean space $\mathbb{R}^k$ for some $k$. We have seen in Proposition \ref{sg_infty}(2) that if $\pi_1(M,x)$ is not finitely generated, then
\[\big(r_i^{-1}\widetilde{M},\tilde{x},\pi_1(M,x)\big)\overset{{\rm GH}}\longrightarrow \big(\mathbb{R}^k,0,G\big)\]
has a limit orbit $G\cdot 0$ that is not connected. If we change the scale $r_i$ to $sr_i$, where $s>0$, then the corresponding limit $\big(s^{-1}\mathbb{R}^k,0,G\big)$ is not equivariantly isometric to $\big(\mathbb{R}^k,0,G\big)$, because the distance between two components of $G\cdot 0$ is changed by a scale of $s^{-1}$. On the other hand, because we assumed that the asymptotic cone of $\widetilde{M}$ is independent of the scaling sequence, one may wonder whether the equivariant limit is independent of the scales as well. If this turns out to be true, then any equivariant asymptotic cone at infinity ought to have a connected orbit at the base point and finite generation would follow. Moreover, it can be shown that if~$G$ is nilpotent and $(s\mathbb{R}^k,0,G)$ are all equivariantly isometric to each other for all scales $s>0$, then the $G=\mathbb{R}^l\times K$, where $K$-action fixes $0$ and the subgroup $\mathbb{R}^l\times \{e\}$ acts as translations. Hence one may expect that if $\widetilde{M}$ has a unique asymptotic cone as $\mathbb{R}^k$, then any equivariant asymptotic cone $\big(\mathbb{R}^k,0,G\big)$ of $\big(\widetilde{M},\pi_1(M)\big)$ should have a~splitting structure $G=\mathbb{R}^l\times K$ as described above, where~$\pi_1(M)$ is nilpotent.

The author studied Question \ref{equi_stab_infty} when the Riemannian universal cover satisfies one of the conditions below.
\begin{itemize}\itemsep=0pt
	\item $\widetilde{M}$ is $k$-Euclidean at infinity, where $k$ is an integer. This means that any asymptotic cone of $\widetilde{M}$ is a metric cone that splits off exactly~$\mathbb{R}^k$.
	\item $\widetilde{M}$ is $(C(X),\epsilon_X)$-stable at infinity, where $C(X)\in\mathcal{M}(n,0)$ is a~metric cone and~$\epsilon_X$ is a~constant only depending on the cross-section~$X$. This means that any asymptotic cone of~$\widetilde{M}$ is a metric cone, whose cross-section is $\epsilon_X$-close to $X$ with respect to Gromov--Hausdorff distance.
\end{itemize}

\begin{Theorem}[\cite{Pan2,Pan3}]\label{cover_fg}
	Let $M$ be an open $n$-manifold of $\operatorname{Ric}\ge 0$. Suppose that one of the following conditions holds on the Riemannian universal cover $\widetilde{M}$:
\begin{enumerate}\itemsep=0pt
\item[$(1)$] $\widetilde{M}$ is $k$-Euclidean at infinity for some $0\le k\le n$; or
\item[$(2)$] $\widetilde{M}$ is $(C(X),\epsilon_X)$-stable at infinity for some metric cone $C(X)\in \mathcal{M}(n,0)$ and a sufficiently small $\epsilon_X>0$.
\end{enumerate}
	Then $\pi_1(M)$ is finitely generated.
\end{Theorem}

Both cases include the case that $M$ has nonnegative sectional curvature, with which $\widetilde{M}$ has a unique asymptotic cone as a metric cone. They also cover the case that $\widetilde{M}$ has Euclidean volume growth and a unique asymptotic cone. Another special case for~(2) is when $\widetilde{M}$ has Euclidean volume growth of constant at least $1-\epsilon(n)$ for some small $\epsilon(n)>0$; in this case, due to Cheeger--Colding almost maximal volume rigidity~\cite{CC1}, any asymptotic cone of $\widetilde{M}$ has a~cross-section being Gromov--Hausdorff close to the unit $(n-1)$-dimensional sphere.

We view the condition in Theorem \ref{cover_fg} as a geometric stability condition of $\widetilde{M}$ at infinity. As indicated, the key to proving finite generation is establishing corresponding equivariant stability at infinity. For simplicity, we only state the result of equivariant stability when $\widetilde{M}$ has a unique cone as a metric cone $C(Z)$ here.

\begin{Theorem}\label{unique_cone_eq}
	Let $M$ be an open $n$-manifold with $\operatorname{Ric}\ge 0$ and a nilpotent fundamental group. Suppose that the universal cover $\widetilde{M}$ has a unique asymptotic cone as a metric cone $C(Z)$. Then there exist a closed nilpotent subgroup $K$ of $\operatorname{Isom}(Z)$ and an integer $l\in[0,n]$ such that any equivariant asymptotic cone $\big(\widetilde{Y},\tilde{y},G\big)$ of $\big(\widetilde{M},\pi_1(M,x)\big)$ is $\big(C(Z),z,\mathbb{R}^l\times K\big)$, where $K$-action fixes $z$ and the subgroup $\mathbb{R}^l\times \{e\}$ acts as translations in the Euclidean factor of $C(Z)$.
\end{Theorem}

To prove Theorem~\ref{unique_cone_eq}, it is essential to treat the set of all equivariant asymptotic cones of $\big(\widetilde{M},\pi_1(M,x)\big)$, denoted as $\Omega\big(\widetilde{M},\pi_1(M,x)\big)$, as a whole, not as individual spaces. It is known that $\Omega\big(\widetilde{M},\pi_1(M,x)\big)$ is compact and connected in the pointed equivariant Gromov--Hausdorff topology. This fact serves as an intuition behind the proof.

The proof of Theorem \ref{unique_cone_eq} is a contradicting argument. Suppose that there are two different equivariant asymptotic cones $(C(Z),z,G_1)$ and $(C(Z),z,G_2)$. Using the connectedness of $\Omega\big(\widetilde{M},\pi_1(M,x)\big)$, we can find a chain of spaces in $\Omega\big(\widetilde{M},\pi_1(M,x)\big)$ connecting $(C(Z),z,G_1)$ and $(C(Z),z,G_2)$, such that the adjacent spaces in the chain are very close in the equivariant Gromov--Hausdorff topology. Roughly speaking, we find a contradiction somewhere in the chain. The author developed two technical tools to achieve this. One is a critical rescaling argument. We view $(C(Z),z,G_1)$ and $(C(Z),z,G_2)$ as the limits coming from two different scales
\[\big(r_i^{-1}\widetilde{M},\tilde{x},\pi_1(M,x)\big)\overset{{\rm GH}}\longrightarrow(C(Z),z,G_1),\qquad \big(s_i^{-1}\widetilde{M},\tilde{x},\pi_1(M,x)\big)\overset{{\rm GH}}\longrightarrow(C(Z),z,G_2).\]
Adjusting the scales, we can assume that one is the rescaling of the other by some $t_i\to\infty$. The critical rescaling argument chooses a suitable intermediate scaling sequence $l_i\to\infty$ with $l_i\le t_i$ and looks for a contradiction in the corresponding limit coming from~$l_i$. This argument helps reduce the problem to a compact one regarding isometric actions on the cross-section $Z$. The other tool is an equivariant Gromov--Hausdorff distance gap between isometric actions on $Z$. To describe this, we think of a chain of spaces $\{(Z,H_j)\}_{j=1}^m$, where each $H_j$ is a closed subgroup of~$\operatorname{Isom}(Z)$. The equivariant Gromov--Hausdorff distance gap implies that if the distance between adjacent spaces is sufficiently small along the chain, then all~$(Z,H_j)$ have to be the same space. Theorem~\ref{isom_Lie} is crucial in establishing such a distance gap.

Based on \cite{Pan2,Pan3}, we ask the following question on open manifolds whose Riemannian universal covers have Euclidean volume growth.

\begin{Question}\label{ques_cover_eu}When the Riemannian universal cover $\widetilde{M}$ has Euclidean volume growth, is there any equivariant stability among all equivariant asymptotic cones of $\big(\widetilde{M},\pi_1(M,x)\big)$?
\end{Question}

One may think of the Euclidean volume growth condition as certain geometric stability condition at infinity as well, but a very weak one: we only know that all asymptotic cones are metric cones whose cross sections share the same $(n-1)$-dimensional Hausdorff volume. An affirmative answer to Question~\ref{ques_cover_eu} would imply the finite generation for this class of manifolds.

\section{Virtual abelianness}\label{sec_abel}

\subsection{Wei's construction and related examples}\label{sec_exmp}

Wei first constructed an open manifold of positive Ricci curvature whose fundamental group is not virtually abelian, as a torsion-free nilpotent group \cite{Wei}. Based on this construction, later Wilking showed that any finitely generated virtually nilpotent group can be realized as the fundamental group of some open manifold with positive Ricci curvature \cite{Wilk}. We briefly describe Wei's construction and other related examples.

Let $\big(S^{p-1},ds^2\big)$ be the unit $(p-1)$-dimensional sphere and let $(N,g_0)$ be a closed manifold. A doubly warped product $[0,\infty)\times_f S^{p-1}\times_h N$ has the metric
\[g=dr^2+f(r)^2ds^2+h(r)^2g_0,\]
It is diffeomorphic to $\mathbb{R}^p\times N$ if the warping functions $f$ and $h$ are smooth and satisfy
\[ f(0)=0,\qquad f'(0)=1,\qquad f''(0)=0,\qquad h(0)>0,\qquad h'(0)=0\]
and $f(r)>0$, $h(r)>0$ for all $r>0$.

Nabonnand constructed a doubly warped product $[0,\infty)\times_f S^{2}\times_h S^1$ of positive Ricci curvature~\cite{Nab}. Later, Bergery generalized this method to a doubly warped product $M=[0,\infty)\times_f S^{p-1}\times_h N$ of positive Ricci curvature, where $N$ is any closed manifold with nonnegative Ricci curvature~\cite{Bergery}. Note that $\pi_1(M)=\pi_1(N)$ is always virtually abelian. In this construction, warping function $h(r)$ can be adjusted so that $h(r)\to 0$ or $h(r)\to c>0$ as~$r\to\infty$.

We illustrate Wei's example. Let $\widetilde{N}$ be a simply connected nilpotent Lie group and let $\Gamma$ be a lattice in $\widetilde{N}$. When $\widetilde{N}$ has nilpotency length $\ge 2$, $\Gamma$ is not virtually abelian and thus the compact quotient manifold $N=\widetilde{N}/\Gamma$ does not admit a metric of nonnegative Ricci curvature. However, $N$ admits a family of metrics $g_r$ with almost nonnegative sectional curvature. Wei constructed a warped product $M=[0,\infty)\times_f S^{p-1}\times N_r$ of positive Ricci curvature with the metric
\[g=dr^2+f(r)^2 ds^2+g_r.\]
This metric is not a doubly warped product in the usual sense since the metric $g_r$ decays at different rates for different steps of~$N$. Note that $\pi_1(M)=\pi_1(N)$ is not virtually abelian. In Wei's construction, $\operatorname{diam}(N,g_r)$ has polynomial decay as $r\to\infty$.

We mention that one can construct a doubly warped product $[0,\infty)\times_f S^{p-1}\times_g S^1$ with positive Ricci curvature and a logarithm decaying $h(r)$ (see Appendix~B of~\cite{Pan4}). However, unlike Wei's construction, it is impossible to use logarithm decaying functions to warp a nilpotent manifold $M=[0,\infty)\times_f S^{p-1}\times N_r$ while keeping positive Ricci curvature. We will see the reason in the next section.

\subsection{Escape rate}\label{sec_escape}

In the light of Theorem \ref{bdd_vir_abel}, we see that the virtual abelianness is related to where the representing geodesic loops of $\pi_1(M,x)$ are positioned. We also mentioned in Section \ref{sec_split} that the escape phenomenon is prevalent for nonnegative Ricci curvature. Motivated by these, the author introduced the \textit{escape rate} to quantify this phenomenon by comparing the size of representing geodesic loops to their length \cite{Pan4}.

\begin{Definition}\label{def_escape}
Let $(M,x)$ be an open manifold with an infinite fundamental group. We define the \textit{escape rate} of $(M,x)$, a scaling invariant, as
\[E(M,x)=\limsup_{|\gamma|\to \infty} \dfrac{d_{\rm H}(x,c_\gamma)}{|\gamma|},\]
	where $\gamma\in\pi_1(M,x)$, $|\gamma|=\operatorname{length}(c_\gamma)$, and $d_{\rm H}$ is the Hausdorff distance.
\end{Definition}

As a loop based at $x$, its size~$d_{\rm H}(x,c_\gamma)$ is at most half of its length. Hence~$E(M,x)$ always ranges from~$0$ to~$1/2$. When all representing geodesic loops stay inside a bounded set or escape at a sublinear rate compared to their length, we have $E(M,x)=0$.

Regarding examples as warped products $[0,\infty)\times_f S^{p-1}\times_h N$, their escape rates are determined by the warping function $h(r)$. As $\operatorname{diam}(N_r)$ decreases to $0$, the representing geodesic loop will take advantage of the thin end to shorten its length, which will in turn enlarge its size. In other words, one should expect that the faster $\operatorname{diam}(N_r)$ decays, the larger the escape rate is. In fact, if $h(r)$ has polynomial decay, as in Wei's construction, then $E(M,x)$ is positive; if $h(r)$ has logarithm decay or converges to a positive constant, then $E(M,x)=0$.

It turns out that this simple geometric quantity can measure the structure of fundamental group. First recall that if $\pi_1(M,x)$ is not finitely generated, then we have infinitely many short generators whose representing geodesic loops are minimal up to halfway. This shows that $E(M,x)=1/2$ if $\pi_1(M)$ is not finitely generated. The author proved that zero escape rate implies virtual abelianness, which substantially generalizes Theorem \ref{bdd_vir_abel}.

\begin{Theorem}[\cite{Pan4}]\label{zero_escape}
	Let $(M,x)$ be an open $n$-manifold of $\operatorname{Ric}\ge 0$. If $E(M,x)=0$, then $\pi_1(M,x)$ is virtually abelian.
\end{Theorem}

We outline the proof of Theorem \ref{zero_escape}, which goes through the asymptotic geometry. A crucial step is showing that the following statements are equivalent:
\begin{enumerate}\itemsep=0pt
\item[(1)] $E(M,x)=0$;
\item[(2)] in any equivariant asymptotic cone $\big(\widetilde{Y},\tilde{y},G\big)$ of $\big(\widetilde{M},\pi_1(M,x)\big)$, the orbit $G\cdot \tilde{y}$ is geodesic in $\widetilde{Y}$, that is, its intrinsic and extrinsic metric agree;
\item[(3)] in any equivariant asymptotic cone $\big(\widetilde{Y},\tilde{y},G\big)$ of $\big(\widetilde{M},\pi_1(M,x)\big)$, the orbit $G\cdot \tilde{y}$ is geodesic in $\widetilde{Y}$ and is isometric to a standard Euclidean space.
\end{enumerate}

The proof of (2) $\Rightarrow$ (3) relies on the Cheeger--Colding splitting theorem and a critical rescaling argument, which is the type of argument we mentioned before in Section~\ref{sec_equi_stab} when explaining Theorem~\ref{unique_cone_eq}. The details are distinct from the proof of Theorem~\ref{unique_cone_eq} due to different contexts.

To derive virtual abelianness from (3), we take a nilpotent subgroup $N$ in $\pi_1(M,x)$ of finite index. For any $r_i\to\infty$, we consider
\[ \big(r_i^{-1}\widetilde{M},\tilde{x},N,\pi_1(M,x)\big)\overset{{\rm GH}}\longrightarrow \big(\widetilde{Y},\tilde{y},H,G\big).\]
It can be shown that $H$ acts on $G\cdot \tilde{y}=H\cdot \tilde{y}=\mathbb{R}^k$ by translations. This limiting behavior restricts $N$-action at large scale on~$\widetilde{M}$: they acts as almost translations in the sense that $d\big(\gamma^2\tilde{x},\tilde{x}\big)$ is close to twice of $d(\gamma\tilde{x},\tilde{x})$ when $\gamma$ has large displacement at $\tilde{x}$. It was observed in~\cite{Pan3} that this feature implies virtual abelianness.

We explain why the Heisenberg $3$-group $H^3$ cannot act this way isometrically. For convenience, we write $|\gamma|=d(\gamma\tilde{x},\tilde{x})$. The commutator calculus in $H^3$ has $\big[g^k,h\big]=[g,h]^k=\big[g,h^k\big]$ for all $g,h\in H^3$ and $k\in\mathbb{N}$. Now let $g,h\in H^3$ with $[g,h]\not= e$. We choose a large integer~$k$ such that $R\le \big|[g,h]^{(k^2)}\big|=\big|\big[g^k,h^k\big]\big|$. We set $\alpha=g^k$ and $\beta=h^k$, then continue to raise the power. For any integer $p$ and $l=2^p$, we have
\[\big|[\alpha,\beta]^{(l^2)}\big|=\big|\big[\alpha^l,\beta^l\big]\big|\le 2\cdot 2^p(|\alpha|+|\beta|).\]
On the other hand, the almost translation at large scale condition implies that
\[\big|[\alpha,\beta]^{(l^2)}\big|\ge (1.9)^{2p} |[\alpha,\beta]|.\]
This leads to a contradiction when $p$ is sufficiently large.

For manifolds in Theorem~\ref{cover_fg}, as a consequence of their equivariant stability at infinity, they have zero escape rate~\cite{Pan4}. We can further bound the index of the abelian subgroup if in addition the universal cover has Euclidean volume growth.

\begin{Theorem}[\cite{Pan3}]\label{cover_abel}
	Let $M$ be an open $n$-manifold of $\operatorname{Ric}\ge 0$. Suppose that one of the following conditions holds on the Riemannian universal cover $\widetilde{M}$:
\begin{enumerate}\itemsep=0pt
\item[$(1)$] $\widetilde{M}$ is $k$-Euclidean at infinity for some $0\le k\le n$; or
\item[$(2)$] $\widetilde{M}$ is $(C(X),\epsilon_X)$-stable at infinity for some metric cone $C(X)\in \mathcal{M}(n,0)$ and a sufficiently small $\epsilon_X>0$.
\end{enumerate}
	Then $\pi_1(M)$ is virtually abelian. If in addition $\widetilde{M}$ has Euclidean volume growth of constant at least $L>0$, then the index can be bounded by some constant $C(n,L)$.
\end{Theorem}

For manifolds whose universal covers have Euclidean volume growth, if one can answer Question \ref{ques_cover_eu} affirmatively, then that would confirm $E(M,x)=0$ and the conjecture below.

\begin{Conjecture}\label{conj_cover_eu_vol}
	Given $n$ and $L\in(0,1]$, there exists a constant $C(n,L)$ such that the following holds.
	
	Let $M$ be an open $n$-manifold of $\operatorname{Ric}\ge 0$. If the Riemannian universal cover of $M$ has Euclidean volume growth of constant at least $L$, then $\pi_1(M)$ is finitely generated and contains an abelian subgroup of index at most $C(n,L)$.
\end{Conjecture}

As mentioned in Section \ref{sec_equi_stab}, the case $L$ being close to $1$ is indeed included in Theorems~\ref{cover_fg}(2) and~\ref{cover_abel}(2).

Theorem~\ref{zero_escape} is further generalized to an escape rate gap by the author.

\begin{Theorem}[\cite{Pan5}]\label{escape_gap}
Given $n$, there exists a constant $\epsilon(n)>0$ such that the following holds.
	
Let $(M,x)$ be an open $n$-manifold of $\operatorname{Ric}\ge 0$. If $E(M,x)\le \epsilon(n)$, then $\pi_1(M,x)$ is virtually abelian.
\end{Theorem}

To end this section, we pose the question whether escape rate can detect nilpotent groups of different steps.

\begin{Question}For each $n$, are there positive constants $\epsilon(n,k)$ strictly increasing in $k$ such that the following holds?
	
For any open $n$-manifold $(M,x)$ of $\operatorname{Ric}\ge 0$, if $E(M,x)\le \epsilon(n,k)$, then $\pi_1(M)$ contains a free nilpotent subgroup of nilpotency length at most $k$ with finite index.
\end{Question}

\subsection*{Acknowledgements}

The author is partially supported by AMS Simons travel grant. The author would like to thank Guofang Wei and Christina Sormani for their encouragement and suggestions when preparing this survey.

\pdfbookmark[1]{References}{ref}
\LastPageEnding


\begin{thebibliography}{99}
\footnotesize\itemsep=0pt

\bibitem{AG}
Abresch U., Gromoll D., On complete manifolds with nonnegative {R}icci
 curvature, \href{https://doi.org/10.2307/1990957}{\textit{J.~Amer. Math. Soc.}} \textbf{3} (1990), 355--374.

\bibitem{An}
Anderson M.T., On the topology of complete manifolds of nonnegative {R}icci
 curvature, \href{https://doi.org/10.1016/0040-9383(90)90024-E}{\textit{Topology}} \textbf{29} (1990), 41--55.

\bibitem{AKL}
Anderson M.T., Kronheimer P.B., LeBrun C., Complete {R}icci-flat {K}\"ahler
 manifolds of infinite topological type, \href{https://doi.org/10.1007/BF01228345}{\textit{Comm. Math. Phys.}}
 \textbf{125} (1989), 637--642.

\bibitem{Bergery}
B\'erard-Bergery L., Quelques exemples de vari\'et\'es riemanniennes
 compl\`etes non compactes \`a courbure de {R}icci positive,
 \textit{C.~R.~Acad. Sci. Paris S\'er.~I Math.} \textbf{302} (1986), 159--161.

\bibitem{Bishop}
Bishop R.L., A relation between volume, mean curvature and diameter,
 \textit{Amer. Math. Soc. Not.} \textbf{10} (1963), 364--364.

\bibitem{CC1}
Cheeger J., Colding T.H., Lower bounds on {R}icci curvature and the almost
 rigidity of warped products, \href{https://doi.org/10.2307/2118589}{\textit{Ann. of Math.}} \textbf{144} (1996),
 189--237.

\bibitem{CC2}
Cheeger J., Colding T.H., On the structure of spaces with {R}icci curvature
 bounded below.~{I}, \href{https://doi.org/10.4310/jdg/1214459974}{\textit{J.~Differential Geom.}} \textbf{46} (1997),
 406--480.

\bibitem{CC3}
Cheeger J., Colding T.H., On the structure of spaces with {R}icci curvature
 bounded below.~{II}, \href{https://doi.org/10.4310/jdg/1214342145}{\textit{J.~Differential Geom.}} \textbf{54} (2000),
 13--35.

\bibitem{CC4}
Cheeger J., Colding T.H., On the structure of spaces with {R}icci curvature
 bounded below.~{III}, \href{https://doi.org/10.4310/jdg/1214342146}{\textit{J.~Differential Geom.}} \textbf{54} (2000),
 37--74.

\bibitem{CG_split}
Cheeger J., Gromoll D., The splitting theorem for manifolds of nonnegative
 {R}icci curvature, \href{https://doi.org/10.4310/jdg/1214430220}{\textit{J.~Differential Geometry}} \textbf{6} (1971),
 119--128.

\bibitem{CG_soul}
Cheeger J., Gromoll D., On the structure of complete manifolds of nonnegative
 curvature, \href{https://doi.org/10.2307/1970819}{\textit{Ann. of Math.}} \textbf{96} (1972), 413--443.

\bibitem{ChN1}
Cheeger J., Naber A., Lower bounds on {R}icci curvature and quantitative
 behavior of singular sets, \href{https://doi.org/10.1007/s00222-012-0394-3}{\textit{Invent. Math.}} \textbf{191} (2013),
 321--339, \href{https://arxiv.org/abs/1103.1819}{arXiv:1103.1819}.

\bibitem{ChN2}
Cheeger J., Naber A., Regularity of {E}instein manifolds and the codimension 4
 conjecture, \href{https://doi.org/10.4007/annals.2015.182.3.5}{\textit{Ann. of Math.}} \textbf{182} (2015), 1093--1165,
 \href{https://arxiv.org/abs/1406.6534}{arXiv:1406.6534}.

\bibitem{Co}
Colding T.H., Ricci curvature and volume convergence, \href{https://doi.org/10.2307/2951841}{\textit{Ann. of Math.}}
 \textbf{145} (1997), 477--501.

\bibitem{CoN}
Colding T.H., Naber A., Sharp {H}\"{o}lder continuity of tangent cones for
 spaces with a lower {R}icci curvature bound and applications, \href{https://doi.org/10.4007/annals.2012.176.2.10}{\textit{Ann. of
 Math.}} \textbf{176} (2012), 1173--1229, \href{https://arxiv.org/abs/1102.5003}{arXiv:1102.5003}.

\bibitem{EM}
Evans B., Moser L., Solvable fundamental groups of compact {$3$}-manifolds,
 \href{https://doi.org/10.2307/1996169}{\textit{Trans. Amer. Math. Soc.}} \textbf{168} (1972), 189--210.

\bibitem{Fky}
Fukaya K., Theory of convergence for {R}iemannian orbifolds, \href{https://doi.org/10.4099/math1924.12.121}{\textit{Japan.~J.
 Math. (N.S.)}} \textbf{12} (1986), 121--160.

\bibitem{FY}
Fukaya K., Yamaguchi T., The fundamental groups of almost non-negatively curved
 manifolds, \href{https://doi.org/10.2307/2946606}{\textit{Ann. of Math.}} \textbf{136} (1992), 253--333.

\bibitem{Gro1}
Gromov M., Almost flat manifolds, \href{https://doi.org/10.4310/jdg/1214434488}{\textit{J.~Differential Geometry}} \textbf{13}
 (1978), 231--241.

\bibitem{Gro2}
Gromov M., Groups of polynomial growth and expanding maps, \href{https://doi.org/10.1007/BF02698687}{\textit{Inst. Hautes
 \'Etudes Sci. Publ. Math.}} (1981), 53--73.

\bibitem{Gro_book}
Gromov M., Metric structures for {R}iemannian and non-{R}iemannian spaces,
 \textit{Modern Birkh\"auser Classics}, \href{https://doi.org/10.1007/978-0-8176-4583-0}{Birkh\"auser Boston, Inc.}, Boston, MA, 2007.

\bibitem{Honda}
Honda S., On low-dimensional {R}icci limit spaces, \href{https://doi.org/10.1017/S0027763000010667}{\textit{Nagoya Math.~J.}}
 \textbf{209} (2013), 1--22.

\bibitem{KW}
Kapovitch V., Wilking B., Structure of fundamental groups of manifolds of
 {R}icci curvature bounded below, \href{https://arxiv.org/abs/1105.5955}{arXiv:1105.5955}.

\bibitem{Li}
Li P., Large time behavior of the heat equation on complete manifolds with
 nonnegative {R}icci curvature, \href{https://doi.org/10.2307/1971385}{\textit{Ann. of Math.}} \textbf{124} (1986),
 1--21.

\bibitem{Liu}
Liu G., 3-manifolds with nonnegative {R}icci curvature, \href{https://doi.org/10.1007/s00222-012-0428-x}{\textit{Invent. Math.}}
 \textbf{193} (2013), 367--375, \href{https://arxiv.org/abs/1108.1888}{arXiv:1108.1888}.

\bibitem{Menguy}
Menguy X., Noncollapsing examples with positive {R}icci curvature and infinite
 topological type, \href{https://doi.org/10.1007/PL00001632}{\textit{Geom. Funct. Anal.}} \textbf{10} (2000), 600--627.

\bibitem{Mil}
Milnor J., A note on curvature and fundamental group, \href{https://doi.org/10.4310/jdg/1214501132}{\textit{J.~Differential
 Geometry}} \textbf{2} (1968), 1--7.

\bibitem{Nab}
Nabonnand P., Sur les vari\'et\'es riemanniennes compl\`etes \`a courbure de
 {R}icci positive, \textit{C.~R.~Acad. Sci. Paris S\'er.~A-B} \textbf{291}
 (1980), A591--A593.

\bibitem{Pan2}
Pan J., Nonnegative {R}icci curvature, stability at infinity and finite
 generation of fundamental groups, \href{https://doi.org/10.2140/gt.2019.23.3203}{\textit{Geom. Topol.}} \textbf{23} (2019),
 3203--3231, \href{https://arxiv.org/abs/1710.05498}{arXiv:1710.05498}.

\bibitem{Pan1}
Pan J., A proof of {M}ilnor conjecture in dimension 3, \href{https://doi.org/10.1515/crelle-2017-0057}{\textit{J.~Reine Angew.
 Math.}} \textbf{758} (2020), 253--260, \href{https://arxiv.org/abs/1703.08143}{arXiv:1703.08143}.

\bibitem{Pan3}
Pan J., Nonnegative {R}icci curvature, almost stability at infinity, and
 structure of fundamental groups, \href{https://arxiv.org/abs/1809.10220}{arXiv:1809.10220}.

\bibitem{Pan4}
Pan J., On the escape rate of geodesic loops in an open manifold with
 nonnegative {R}icci curvature, \textit{Geom. Topol.}, {t}o appear,
 \href{https://arxiv.org/abs/2003.01326}{arXiv:2003.01326}.

\bibitem{Pan5}
Pan J., Nonnegative {R}icci curvature and escape rate gap, {i}n preparation.

\bibitem{PR}
Pan J., Rong X., Ricci curvature and isometric actions with scaling
 nonvanishing property, \href{https://arxiv.org/abs/1808.02329}{arXiv:1808.02329}.

\bibitem{Per}
Perelman G., Manifolds of positive {R}icci curvature with almost maximal
 volume, \href{https://doi.org/10.2307/2152760}{\textit{J.~Amer. Math. Soc.}} \textbf{7} (1994), 299--305.

\bibitem{Per_ex}
Perelman G., Construction of manifolds of positive {R}icci curvature with big
 volume and large {B}etti numbers, in Comparison Geometry ({B}erkeley, {CA},
 1993--94), \textit{Math. Sci. Res. Inst. Publ.}, Vol.~30, Cambridge
 University Press, Cambridge, 1997, 157--163.

\bibitem{SchYau}
Schoen R., Yau S.T., Complete three-dimensional manifolds with positive {R}icci
 curvature and scalar curvature, in Seminar on {D}ifferential {G}eometry,
 \textit{Ann. of Math. Stud.}, Vol.~102, Princeton University Press,
 Princeton, N.J., 1982, 209--228.

\bibitem{ShaYang1}
Sha J.-P., Yang D., Examples of manifolds of positive {R}icci curvature,
 \href{https://doi.org/10.4310/jdg/1214442635}{\textit{J.~Differential Geometry}} \textbf{29} (1989), 95--103.

\bibitem{ShaYang2}
Sha J.P., Yang D., Positive {R}icci curvature on the connected sums of
 {$S^n\times S^m$}, \href{https://doi.org/10.4310/jdg/1214446032}{\textit{J.~Differential Geom.}} \textbf{33} (1991),
 127--137.

\bibitem{Shara}
Sharafutdinov V.A., Convex sets in a manifold of nonnegative curvature,
 \href{https://doi.org/10.1007/BF01140282}{\textit{Math. Notes}} \textbf{26} (1979), 556--560.

\bibitem{SS}
Shen Z., Sormani C., The topology of open manifolds with nonnegative {R}icci
 curvature, \textit{Commun. Math. Anal.} (2008), 20--34,
 \href{https://arxiv.org/abs/math.DG/0606774}{arXiv:math.DG/0606774}.

\bibitem{Sor2}
Sormani C., The almost rigidity of manifolds with lower bounds on {R}icci
 curvature and minimal volume growth, \href{https://doi.org/10.4310/CAG.2000.v8.n1.a6}{\textit{Comm. Anal. Geom.}} \textbf{8}
 (2000), 159--212, \href{https://arxiv.org/abs/math.DG/9903171}{arXiv:math.DG/9903171}.

\bibitem{Sor}
Sormani C., Nonnegative {R}icci curvature, small linear diameter growth and
 finite generation of fundamental groups, \href{https://doi.org/10.4310/jdg/1214339792}{\textit{J.~Differential Geom.}}
 \textbf{54} (2000), 547--559, \href{https://arxiv.org/abs/math.DG/9809133}{arXiv:math.DG/9809133}.

\bibitem{Sor3}
Sormani C., On loops representing elements of the fundamental group of a
 complete manifold with nonnegative {R}icci curvature, \href{https://doi.org/10.1512/iumj.2001.50.2048}{\textit{Indiana Univ.
 Math.~J.}} \textbf{50} (2001), 1867--1883, \href{https://arxiv.org/abs/math.DG/9904096}{arXiv:math.DG/9904096}.

\bibitem{SW}
Sormani C., Wei G., Various covering spectra for complete metric spaces,
 \href{https://doi.org/10.4310/AJM.2015.v19.n1.a7}{\textit{Asian~J. Math.}} \textbf{19} (2015), 171--202, \href{https://arxiv.org/abs/1211.7123}{arXiv:1211.7123}.

\bibitem{Wan}
Wan J., On the fundamental group of complete manifolds with almost {E}uclidean
 volume growth, \href{https://doi.org/10.1090/proc/14571}{\textit{Proc. Amer. Math. Soc.}} \textbf{147} (2019),
 4493--4498, \href{https://arxiv.org/abs/1902.05292}{arXiv:1902.05292}.

\bibitem{Wei}
Wei G., Examples of complete manifolds of positive {R}icci curvature with
 nilpotent isometry groups, \href{https://doi.org/10.1090/S0273-0979-1988-15653-4}{\textit{Bull. Amer. Math. Soc. (N.S.)}} \textbf{19}
 (1988), 311--313.

\bibitem{Wilk}
Wilking B., On fundamental groups of manifolds of nonnegative curvature,
 \href{https://doi.org/10.1016/S0926-2245(00)00030-9}{\textit{Differential Geom. Appl.}} \textbf{13} (2000), 129--165.

\bibitem{Wraith1}
Wraith D., Surgery on {R}icci positive manifolds, \href{https://doi.org/10.1515/crll.1998.082}{\textit{J.~Reine Angew.
 Math.}} \textbf{501} (1998), 99--113.

\bibitem{Wraith2}
Wraith D., New connected sums with positive {R}icci curvature, \href{https://doi.org/10.1007/s10455-007-9066-8}{\textit{Ann.
 Global Anal. Geom.}} \textbf{32} (2007), 343--360.

\bibitem{Wu}
Wu B.Y., On the fundamental group of {R}iemannian manifolds with nonnegative
 {R}icci curvature, \href{https://doi.org/10.1007/s10711-012-9730-4}{\textit{Geom. Dedicata}} \textbf{162} (2013), 337--344.

\bibitem{Wylie}
Wylie W.C., Noncompact manifolds with nonnegative {R}icci curvature,
 \href{https://doi.org/10.1007/BF02922066}{\textit{J.~Geom. Anal.}} \textbf{16} (2006), 535--550,
 \href{https://arxiv.org/abs/math.DG/0510139}{arXiv:math.DG/0510139}.

\end{thebibliography}
\end{document}